\documentclass{amsart}
\usepackage{amsmath,amssymb}
\usepackage{stmaryrd}
\numberwithin{equation}{section}
\theoremstyle{plain}
\newtheorem{theorem}{Theorem}[section]
\newtheorem{lemma}[theorem]{Lemma}
\newtheorem{proposition}[theorem]{Proposition}
\newtheorem{corollary}[theorem]{Corollary}
\theoremstyle{definition}
\newtheorem{definition}[theorem]{Definition}

\newtheorem{remark}[theorem]{Remark}
\newtheorem{example}[theorem]{Example}
\renewcommand{\l}{\mathcal{L}}
\renewcommand{\t}{\times}

\newcommand{\lp}{\left(}
\newcommand{\rp}{\right)}
\newcommand{\lac}{\left\{}
\newcommand{\rac}{\right\}}
\newcommand{\lcr}{\left[}
\newcommand{\rcr}{\right]}
\newcommand{\abs}[1]{\lvert#1\rvert}
\newcommand{\s}{\Sigma}
\newcommand{\Sep}{\mathcal{S}}
\newcommand{\sep}{\mathsf{S}}
\newcommand{\mo}{\mathsf{MO}}
\newcommand{\h}{\mathcal H}
\newcommand{\pro}{\mathsf{P}}
\newcommand{\prs}{\mathrel{\s_{_\Downarrow}}}
\newcommand{\p}{\perp}
\newcommand{\pr}{\pi^{-1}}
\newcommand{\ots}{\mathrel{\Downarrow\kern-1.08em\bigcirc}}
\newcommand{\otssub}{{\mathrel{\Downarrow\kern-.63em\bigcirc}}}

\DeclareMathOperator{\aut}{\mathsf{Aut}}
\DeclareMathOperator{\uni}{\mathsf{U}}
\begin{document}
\title{Orthocomplemented weak tensor products}
\author{Boris Ischi}
\email{boris.ischi@edu.ge.ch}
\address{Coll\`ege de Candolle\\ 5 rue d'Italie\\ 1204 Geneva\\ Switzerland}
\begin{abstract}
Let $\l_1$ and $\l_2$ be complete atomistic lattices. In a previous 
paper, we have defined a set $\sep=\sep(\l_1,\l_2)$ of complete 
atomistic lattices, the elements of which are called weak tensor products 
of $\l_1$ and $\l_2$. $\sep$ is defined by means of three axioms, natural 
regarding the description of some compound systems 
in quantum logic. It has been proved that $\sep$ is a complete lattice. The
top element of $\sep$, denoted by $\l_1\ovee\l_2$, is the tensor product of Fraser whereas
the bottom element, denoted by $\l_1\owedge\l_2$, is the {\it box product} of Gr\"atzer 
and Wehrung. With some additional hypotheses on $\l_1$ and $\l_2$ (true for 
instance if $\l_1$ and $\l_2$ are moreover orthomodular with the covering property)
we prove that $\sep$ is a singleton
if and only if $\l_1$ or $\l_2$ is distributive, if and only if
$\l_1\ovee\l_2$ has the covering property. Our main result reads:
$\l\in\sep$ admits an orthocomplementation if and only if
$\l=\l_1\owedge\l_2$. 
At the end, we construct an example
$\l_1\ots\l_2$ in $\sep$ which has the covering property. 
\end{abstract}
\maketitle
\section{Introduction}

Let $\s$ be a non-empty set. By a {\it simple closure space} $\l$
on $\s$, we mean a set of subsets of $\s$, ordered by
set-inclusion, closed under arbitrary set-intersections ({\it
i.e.}, for all $\omega\subseteq \l$, $\bigcap\omega\in\l$), and
containing $\s$, $\emptyset$, and all singletons. We denote the
bottom ($\emptyset$) and top ($\s$) elements by $0$ and $1$
respectively. For $p\in\s$, we identify $p$ with $\{p\}\in\l$.
Hence $p\bigcup q$ stands for $\{p,q\}$.

Let $\l$ be a simple closure space on a (nonempty) set $\s$.
Then $\l$ is a complete atomistic lattice, the atoms of which
correspond to the points ({\it i.e.}, singletons) of $\s$. Note
that if $A\subseteq\s$, then
$\bigvee_\l(A)=\bigcap\{b\in\l\ \lvert\ A\subseteq b\}$.
Conversely, let $\l$ be a complete atomistic lattice. Let $\s$
denote the set of atoms of $\l$, and, for each $a\in\l$, let
$\s[a]$ denote the set of atoms under $a$. Then
$\{\s[a]\ \lvert\ a\in\l\}$ is a simple closure space on $\s$,
isomorphic to $\l$.
For simplicity, we shall from now on deal only with simple closure spaces instead 
of complete atomistic lattices.

Let $\l_1$, $\l_2$ and $\l$ be simple closure spaces on $\s_1$, $\s_2$ and $\s$ respectively.
Then, $\l$ is a {\it weak tensor product} of $\l_1$ and $\l_2$ if
\begin{enumerate}
\item[(P1)] $\s=\s_1\t\s_2$
\item[(P2)] $a_1\t\s_2\cup\s_1\t a_2\in\l$, $\forall$ $a_i\in\l_i$
\item[(P3)] for all $p_1\in\s_1$ and $A_2\subseteq \s_2$, $p_1\t A_2\in\l$ implies $A_2\in\l_2$

\noindent for all $p_2\in\s_2$ and $A_1\subseteq \s_1$, $A_1\t p_2\in\l$ implies $A_1\in\l_1$
\end{enumerate}

Let $\sep(\l_1,\l_2)$ denote the set of weak tensor products of $\l_1$ and $\l_2$, ordered by set-inclusion.
Write $\pi_i:\s_1\t\s_2\rightarrow\s_i$ the projection map defined by $\pi_i(p_1,p_2)=p_i$.
Then, it is easy to check that
\[\l_1\owedge\l_2=\lac\vbox{\vspace{0.5cm}}\bigcap\omega\ \lvert\ \omega\subseteq\lac a_1\t\s_2\cup\s_1\t a_2 \ \lvert\  a_i\in\l_i\rac
,\,\omega\ne\emptyset\rac
\]
is the bottom element of $\sep(\l_1,\l_2)$ whereas
\[\begin{aligned}
\l_1\ovee\l_2=\lac A\subseteq\s_1\t\s_2\ \lvert\ \vbox{\vspace{0.5cm}}\right.&\left.\pi_2\lp p_1\t\s_2\cap A\rp\in\l_2\right.\ \mbox{and}\\
\left.\right.&\left.\pi_1\lp\s_1\times p_2\cap A\rp \in\l_1,\ \forall\ p_i\in\s_i\vbox{\vspace{0.5cm}}\rac
\end{aligned}\]
is the top element of $\sep(\l_1,\l_2)$ \cite{Ischi:2007}.

A cornerstone in quantum logic is the following theorem, by which the Hilbert space
structure of quantum mechanics can be recovered from a certain
number of axioms on the lattice $\l_S$ of experimental propositions concerning
a physical system $S$ \cite{Piron:1964} (see \cite{Maeda/Maeda:handbook}, Theorems 34.5 and 34.9). 
Let $E$ be a
vector space over a $^*-$division ring $\mathbb{K}$. A {\it
Hermitian form} $f$ is a mapping $f$ of $E\t E$ onto $\mathbb{K}$,
linear in the left variable ({\it i.e.}, $f(\lambda_1
x_1+\lambda_2 x_2,y)=\lambda_1 f(x_1,y)+\lambda_2 f(x_2,y)$, where
$x_i, y\in E$ and $\lambda_i\in \mathbb{K}$), such that
$f(x,y)^*=f(y,x)$ and $f(x,x)=0$ implies $x=0$. Let $V$ be a
subspace of $E$. Put $V^\p=\{y\in E\ \lvert\  f(x,y)=0\ \forall x\in
V\}$ and say that $V$ is $E-closed$ if $V=V^{\p\p}$.

\begin{theorem}\label{TheoremPiron}
If $\l_S$ is an irreducible orthocomplemented simple closure space
with the covering property, and of length $\geq 4$, then
there is a $^*-$division ring $\mathbb{K}$ and a vector space
$E$ over $\mathbb{K}$ with an Hermitian form such that $\l_S$ is
ortho-isomorphic to the lattice $\pro(E)$ of $E-$closed
subspaces of $E$. Moreover, if $\mathbb{K}=\mathbb{R}$ or
$\mathbb{C}$ with the usual involution, then $E$ is a Hilbert
space if and only if $\l$ is moreover orthomodular.
\end{theorem}

For some compound quantum systems, called {\it separated systems}, 
it can be shown that $\l_S$, which we
henceforth denote by $\l_{S_{sep}}$, must satisfy
the three axioms P1-P3 defining weak tensor products (see \cite{Ischi:2007}
and references herein) and cannot be isomorphic to the
lattice of closed subspaces of a Hilbert space \cite{AertsHPA84}.
As consequence, some hypotheses of Theorem \ref{TheoremPiron} fail
in $\l_{S_{sep}}$.

In \cite{Aerts:1982}, Aerts proposed 
$\l_1\owedge\l_2$ as a model for $\l_{S_{sep}}$. For $\l_1$ and $\l_2$
orthocomplemented simple closure spaces (in which case, $\l_1\owedge\l_2$
is an orthocomplemented simple closure space), Aerts
proved that if $\l_1\owedge\l_2$ has the covering property
or is orthomodular, then $\l_1$ or $\l_2$ is a power set \cite{Aerts:1982}. As a
consequence, according to Aerts, the covering property and
orthomodularity do not hold in $\l_{S_{sep}}$. A similar
conclusion concerning orthomodularity was obtained by Pulmannov\'a
in \cite{Pulmannova:1985}.

Here we argue that no reasonable model for $\l_{S_{sep}}$ admits
an orthocomplementation.  On the other hand, in case $\l_1$ and $\l_2$ satisfy the axioms 
of theorem \ref{TheoremPiron}, we provide a natural
model for $\l_{S_{sep}}$ which has the covering property. We proceed as follows: following
\cite{Piron:handbook} and \cite{Aerts:1982}, we assume that
$\l_1$, $\l_2$ and $\l_{S_{sep}}$ are simple closure
spaces. 
Let $\l\in\sep(\l_1,\l_2)$.
We prove that if $\l_1$ and $\l_2$ are orthocomplemented
simple closure spaces with the covering property, then $\l$ admits
and orthocomplementation if and only if
$\l=\l_1\owedge\l_2$. We conclude by a simple physical
argument given in \cite{IschiDenver:2004} which shows that
certainly $\l_1\owedge\l_2\subsetneqq\l_{S_{sep}}$.

The rest of the paper is organized as follows. In Section
\ref{SectionMainDefinitions}, we fix some basic terminology and
notation. We define a set $\sep\equiv\sep(\l_1,\cdots,\l_n)$ of
$n$-fold weak tensor products. 

Further, with some additional hypotheses on each $\l_i$ (true for instance 
if
$\l_i$ is moreover orthocomplemented with the covering property) we
prove that $\sep$ is a singleton ({\it i.e.},
$\owedge_i\l_i=\ovee_i\l_i$) if and only if at most one $\l_i$ is
not a power set (Section
\ref{SectionSufficientandNecessaryConditionsFor...}) if and only
if $\owedge_i\l_i$ or $\ovee_i\l_i$ has the covering property
(Section \ref{SectionCoevringProperty}). Note that for the last condition concerning $\ovee_i\l_i$,
we need to make an additional assumption on each $\l_i$.

Finally, Section \ref{SectionOrthocomplementation} is devoted to
our main result and Section \ref{SectionAnExample} to the example
mentioned above.
\section{Main definitions}\label{SectionMainDefinitions}

In this section we give our main definitions. We start with some
background material and basic notations used in the sequel. Parts of this section
are taken directly from \cite{Ischi:2007}.

\begin{definition}
\begin{itemize}
\item
A lattice $\l$ with $0$ and $1$ is {\it orthocomplemented} if
there is a unary operation $\,^\p$ (orthocomplementation), also
denoted by $'$, such that for all $a,\,b\in\l$, $(a^\p)^\p=a$,
$a\leq b$ implies $b^\p\leq a^\p$, and $a\vee a^\p=1$.

\item
An orthocomplemented lattice is {\it orthomodular} if for all
$a,\,b\in\l$, $a\leq b$ implies $b=a\vee(b\land a^\p)$.

\item
A lattice with $0$ has the {\it covering property} if for any atom
$p$ and any $a\in\l$, $p\land a=0$ implies that $p\vee a$
covers $a$ (in symbols $p\vee a\gtrdot a$).

\item
A lattice $\l$ with $0$ and $1$ is called a {\it
DAC-lattice} if $\l$ and its dual $\l^*$ (defined by the converse
order-relation) are atomistic with the covering property
\cite{Maeda/Maeda:handbook}. We say that a lattice $\l$ with $1$
is {\it coatomistic} if the dual $\l^*$  is atomistic.

\item
$2$
denotes the simple closure space isomorphic to the two-element
lattice.

\item
Let $\l$ and $\l_1$ be simple closure spaces on $\s$ and $\s_1$ 
respectively. We write $\aut(\l)$ for the group of
automorphisms of $\l$. Note that any map
$u\colon\l\to\l_1$ sending atoms to atoms induces a mapping
from $\s$ to $\s_1$, which we also call $u$. Thus, if
$u\in\aut(\l)$, then for all $a\in\l$, $u(a)=\{u(p)\ \lvert\ p\in a\}$.

\item
If $\l$ is orthocomplemented, for $p,\, q\in\s$, we write $p\p q$
if $p\in q^\p$, where $q^\p$ stands for $\{q\}^\p$.

\item
If $\h$ is a complex Hilbert space, then $\s_\h$ denotes
the set of $1$-dimensional subspaces of $\h$ and $\pro(\h)$ stands
for the simple closure space on $\s_\h$ isomorphic to the lattice of closed
subspaces of $\h$. Moreover, we write $\uni(\h)$ for the set of
automorphisms of $\pro(\h)$ induced by unitary maps on $\h$.

\item
Finally, we say that $\l$ (respectively $T\subseteq \aut(\l)$)
is {\it transitive} if the action of $\aut(\l)$ (respectively the
action of $T$) on $\s$ is transitive.
\end{itemize}
\end{definition}

\begin{remark}
Note that an orthocomplemented atomistic lattice with the covering
property is a DAC-lattice. Note also that in Theorem
\ref{TheoremPiron}, if instead of orthocomplemented, the simple
closure space is a DAC-lattice, then there is a pair of dual
vector spaces such that a representation theorem similar to
Theorem \ref{TheoremPiron} holds (see \cite{Maeda/Maeda:handbook},
Theorem 33.7).
\end{remark}

\begin{definition}\label{DefinitionAlex}
Let $\{\s_\alpha\}_{\alpha\in\Omega}$ be a
family of nonempty sets, $\mathbf{\s}=\prod_\alpha\s_\alpha$,
$\beta\in\Omega$, $p\in\mathbf{\s}$, $R\subseteq\mathbf{\s}$,
$A\in\prod_\alpha 2^{\s_\alpha}$, and $B\subseteq \s_\beta$. We
shall make use of the following notations:
\begin{enumerate}
\item We denote by $\pi_\beta\colon\mathbf{\s}\to\s_\beta$ the
$\beta-$th coordinate map, {\it i.e.}, $\pi_\beta(p)=p_\beta$.
\item We denote by $p[-,\beta]\colon\s_\beta\to\mathbf{\s}$ the
map that sends $q\in\s_\beta$ to the element of $\mathbf{\s}$
obtained by replacing $p$'s $\beta-$th entry by $q$.
\item We define $R_\beta[p]=\pi_\beta(p[\s_\beta,\beta]\cap
R)$. Note that $R_\beta[p]=\{q\in\s_\beta\ \lvert\  p[q,\beta]\in R\}$.
\item We define $A[B,\beta]\in\prod_\alpha\l_\alpha$ as
$A[B,\beta]_\beta=B$ and $A[B,\beta]_\alpha=A_\alpha$ for
$\alpha\ne\beta$.
\item We write $\overline{A}:=\prod_\alpha A_\alpha$ and
$\overline{A}[B,\beta]:=\overline{A[B,\beta]}$.
\end{enumerate}
We omit the $\beta$ in $p[-,\beta]$ when no confusion can occur.
For instance, we write $p[\s_\beta]$ instead of
$p[\s_\beta,\beta]$.
\end{definition}

\begin{remark}\label{Remarkp[R[p]]=p[S]capR}
$p[R_\beta[p]]=p[\s_\beta]\cap R$.
\end{remark}

\begin{definition}\label{DefinitionPTensorProduct}
Let $\{\l_\alpha\}_{\alpha\in\Omega}$ be a family of simple
closure spaces on $\s_\alpha$. We denote by
$\sep(\l_\alpha,\alpha\in\Omega)$ the set all simple closure spaces
$\l$ on $\mathbf{\s}$ such that
\begin{enumerate}
\item[(P1)] $\mathbf{\s}=\prod_\alpha\s_\alpha$,
\item[(P2)] $\bigcup_\alpha \pr_\alpha(a_\alpha)\in\l$, for all
$a\in\prod_\alpha\l_\alpha$,
\item[(P3)] for all $p\in\mathbf{\s}$, $\beta\in\Omega$, and
$B\subseteq\s_\beta$, we have $p[B,\beta]\in\l$ implies $B\in\l_\beta$.
\end{enumerate}
Let $T=\prod_\alpha T_\alpha$ with $T_\alpha\subseteq
\aut(\l_\alpha)$. We denote by $\Sep_T(\l_\alpha,\alpha\in\Omega)$
the set of all $\l\in\sep(\l_\alpha,\alpha\in\Omega)$ such that
\begin{enumerate}
\item[(P4)] for all $ v\in T$, there is $u\in \aut(\l)$ such that
$u(p)_\alpha=v_\alpha(p_\alpha)$ for all $p\in\mathbf{\s}$ and all
$\alpha\in\Omega$.
\end{enumerate}
We call elements of $\sep(\l_\alpha,\alpha\in\Omega)$ {\it weak
tensor products}.
\end{definition}

\begin{lemma}\label{LemmaAlex}
Let $\{\l_\alpha\}_{\alpha\in\Omega}$ be a family of simple
closure spaces on $\s_\alpha$, $\beta\in\Omega$, and
$\l\in\sep(\l_\alpha,\alpha\in\Omega)$.
\begin{enumerate}
\item For any $a\in\prod_\alpha\l_\alpha$, $\overline{a}\in\l$.
\item For any $b\in\l_\beta$ and $p\in
\prod_\alpha\s_\alpha$, we have $p[b,\beta]\in\l$.
\item For any $B\subseteq\l_\beta$ and $a\in
\prod_\alpha\l_\alpha$, we have $\overline{a}[\bigvee
B,\beta]=\bigvee_{b\in B} \overline{a}[b,\beta]$.
\end{enumerate}
\end{lemma}

\begin{proof} see \cite{Ischi:2007}
\end{proof}

\begin{definition}\label{DefinitionSeparatedProduct}
Let $\{\s_\alpha\}_{\alpha\in\Omega}$ be a family of nonempty sets
and $\{\l_\alpha\subseteq 2^{\s_\alpha}\}_{\alpha\in\Omega}$. Let
$\mathbf{\s}=\prod_\alpha\s_\alpha$. We define
\[\begin{aligned}\mathop{\owedge}_{\alpha\in\Omega}\l_\alpha&:=\lac\vbox{\vspace{0.5cm}}\bigcap \omega\ \lvert\ 
\omega\subseteq \lac\bigcup_\alpha \pr_\alpha(a_\alpha)
\ \lvert\ a\in\prod_\alpha\l_\alpha\rac,\,\omega\ne\emptyset\rac\, ,\\
\mathop{\ovee}_{\alpha\in\Omega}\l_\alpha&:=\lac\vbox{\vspace{0.5cm}} R\subseteq\mathbf{\s}\ \lvert\ 
R_\beta[p]\in\l_\beta,\, \mbox{for all }
p\in\mathbf{\s},\,\beta\in\Omega\rac\, ,
\end{aligned}\]
ordered by set-inclusion.
\end{definition}

\begin{remark} Note that 
$2\owedge\l\cong\l\cong 2\ovee\l$ and $\abs{\sep(2,\l)}=1$.
\end{remark}

We end this section by recalling a definition and two results presented in \cite{Ischi:2007}
that we will use later.

\begin{lemma}\label{LemmaJoinP1diffQ1andP2diffQ2}
Let $\{\l_\alpha\}_{\alpha\in\Omega}$ be a family of simple
closure spaces on $\s_\alpha$ and $\l$ a simple closure space on 
$\mathbf{\s}=\prod_\alpha\s_\alpha$. Suppose that Axiom P2 holds
in $\l$. Let $p,\,q\in\mathbf{\s}$.
\begin{enumerate}
\item If $p_\beta\ne q_\beta$ for at least two $\beta\in\Omega$,
then $p\vee q=p\cup q$.
\item For all $\beta\ne \gamma\in\Omega$ and for all
$b\in\l_\beta$ and $c\in\l_\gamma$ such that $p_\beta\in b$ and
$p_\gamma\in c$, $p[b,\beta]\vee p[c,\gamma]=p[b,\beta]\cup
p[c,\gamma]$.
\end{enumerate}
\end{lemma}

\begin{definition}\label{DefinitionConnected}
Let $\l$ be a simple closure space on $\s$. We say that $\l$ is
{\it weakly connected} if $\l\ne 2$ and if there is a {\it
connected covering} of $\s$, that is a family of subsets
$\{A^\gamma\subseteq\s\ \lvert\ \gamma\in\sigma\}$ such that
\begin{enumerate}
\item $\s=\bigcup \{A^\gamma\ \lvert\ \gamma\in \sigma\}$ and $\abs{
A^\gamma}\geq 2$ for all $\gamma\in\sigma$,
\item for all $\gamma\in\sigma$ and all $p\ne q\in A^\gamma$,
$p\vee q$ contains a third atom,
\item for all $p,\,q\in\s$, there is a finite subset
$\{\gamma_1,\cdots,\gamma_n\}\subseteq\sigma$ such that $p\in
A^{\gamma_1}$, $q\in A^{\gamma_n}$, and such that
$\abs{A^{\gamma_i}\cap A^{\gamma_{i+1}}} \geq 2$ for all $1\leq
i\leq n-1$.
\end{enumerate}
\end{definition}

\begin{remark}
Note that 
weakly connected implies irreducible (see \cite{Ischi:2007}). Finally, let $\l$ be a
simple closure space. Then, if $\l\ne 2$ and $\l$ is irreducible
orthocomplemented with the covering property or an irreducible
DAC-lattice, then $\l$ is weakly connected.
\end{remark}

\begin{theorem}\label{TheoremAutoFactorwithWeaklyConnected}
Let $\Omega$  be a finite set and $\{\l_i\}_{i\in\Omega}$ a finite
family of weakly connected simple closure spaces on $\s_i$. Let
$\l\in\sep(\l_i,i\in\Omega)$ and $u\in\aut(\l)$. Then there is a bijection $f$ of $\Omega$, and for
each $i\in\Omega$, there is an isomorphism $v_i\colon\l_i\to\l_{f(i)}$ 
such that $u(p)_{f(i)}=v_i(p_i)$ for all
$p\in\mathbf{\s}$ and $i\in\Omega$.
\end{theorem}

\begin{remark}\label{remarkread}
It can be useful to have in mind the following pictures in order to
read the proofs below easily. If $\Omega=\{1,2,3\}$, then
\begin{enumerate}
\item $\pr_1(a_1)=a_1\t\s_2\t\s_3$ which we can denote by the symbol $a_1\vert\vert$,
\item $\bigcup_{i=1}^3\pr_i(a_i)=a_1\t\s_2\t\s_3\cup\s_1\t a_2\t\s_3\cup\s_1\t\s_2\t a_3$ which we can denote as
$a_1\vert\vert\cup \vert a_2\vert\cup\vert\vert a_3$. Hence, we can write 
\[
\begin{aligned}
\lp\bigcup_{i=1}^3\pr_i(a_i)\rp\cap\lp\bigcup_{i=1}^3\pr_i(b_i)\rp
=\lcr\vbox{\vspace{0.5cm}} a_1\vert\vert\cup \vert a_2\vert\cup\vert\vert a_3\rcr
\cap\lcr\vbox{\vspace{0.5cm}} b_1\vert\vert\cup \vert b_2\vert\cup\vert\vert b_3\rcr\\
=a_1\cap b_1\vert\vert\cup a_1b_2\vert\cup a_1\vert b_3
\cup b_1a_2\vert\cup \vert a_2\cap b_2\vert\cup \vert a_2b_3
\cup b_1\vert a_3\cup \vert b_2a_3\cup\vert\vert a_3\cap b_3\\
=\bigcup_{f\in3^{\{a,b\}}}\lp\bigcap f^{-1}(1)_1\rp\t\lp\bigcap f^{-1}(2)_2\rp\t\lp\bigcap f^{-1}(3)_3\rp
\end{aligned}
\]
where $a,\,b\in\l_1\t\l_2\t\l_3$.
\end{enumerate}
\end{remark}
\section{Sufficient and necessary conditions for
$\l_1\owedge\l_2=\l_1\ovee\l_2$}\label{SectionSufficientandNecessaryConditionsFor...}

Let $\l_1,\cdots,\l_n$ be simple closure spaces and DAC-lattices. In this section,
we prove that $\owedge_{i=1}^n\l_i=\ovee_{i=1}^n\l_i$ ({\it i.e.},
$\sep(\l_1,\cdots,\l_n)$ has only one element) if and only if
there is $k$ between $1$ and $n$ such that for all $i\ne k$,
$\l_i$ is a power set.

\begin{definition}
Let $\{\l_\alpha\}_{\alpha\in\Omega}$ be a family of simple
closure spaces on $\s_\alpha$ and $\beta\in\Omega$. Let
$\mathbf{\s}=\prod_\alpha\s_\alpha$ and $R\subseteq\mathbf{\s}$.
Then, we define $\bigvee_\beta R:=\bigcup_{p\in\mathbf{\s}}p[\bigvee
R_\beta[p]]$,
where the join is being taken in $\l_\beta$.
\end{definition}

\begin{lemma}\label{LemmaJoininTensorF}
Let $\{\l_\alpha\}_{\alpha\in\Omega}$ be a family of simple
closure spaces on $\s_\alpha$ and
$\l\in\sep(\l_\alpha,\alpha\in\Omega)$. Let
$\mathbf{\s}=\prod_\alpha\s_\alpha$ and $R\subseteq\mathbf{\s}$.
Denote by $\bigvee_\l$ the join in $\l$.
\begin{enumerate}
\item For any $\beta\in\Omega$, $R\subseteq\bigvee_\beta R$ and
$\bigvee_\beta(\bigvee_\beta R)=\bigvee_\beta R$. Moreover, if
$R\subseteq S\subseteq \mathbf{\s}$, then $\bigvee_\beta
R\subseteq\bigvee_\beta S$.
\item For any $f\colon\mathbb{N}\to\Omega$ and any
$n\in\mathbb{N}$, $R^n:=\bigvee_{f(n)}(\cdots(\bigvee_{f(1)}R)\cdots
)\subseteq\bigvee_\l R$.
\item $R\in\ovee_\alpha\l_\alpha$ if and only if $\bigvee_\beta
R=R$ for all $\beta\in\Omega$.
\end{enumerate}
\end{lemma}

\begin{proof} (1) Let $\beta\in\Omega$. First note that
\[
R%
=R\cap\lp\bigcup_{p\in\mathbf{\s}}p[\s_\beta]\rp
=\bigcup_{p\in\mathbf{\s}}(p[\s_\beta]\cap R)%
=\bigcup_{p\in\mathbf{\s}}p[ R_\beta[p]]\, .
\]
As a consequence, since $R_\beta[p]\subseteq\bigvee R_\beta[p]$ in
$\l_\beta$, we have
\[ R=\bigcup_{p\in\mathbf{\s}}p[ R_\beta[p]]
\subseteq \bigcup_{p\in\mathbf{\s}}p[\bigvee
R_\beta[p]]=\bigvee_\beta R,\]
hence, $R\subseteq\bigvee_\beta R$.

Further,
\[
\begin{aligned} \bigvee_\beta(\bigvee_\beta
R)=\bigvee_\beta\lp\bigcup_{p\in\mathbf{\s}}p[\bigvee
R_\beta[p]]\rp=\bigcup_{q\in\mathbf{\s}}q\lcr\bigvee\lp\bigcup_{p\in\mathbf{\s}}p[\bigvee
R_\beta[p]]\rp_\beta[q]\rcr\hfill\\[4mm]
\hfill=\bigcup_{q\in\mathbf{\s}} q\lcr\bigvee\lp q[\bigvee
R_\beta[q]]\rp_\beta[q]\rcr=\bigcup_{q\in\mathbf{\s}}
q\lcr\bigvee\lp \bigvee
R_\beta[q]\rp\rcr=\bigcup_{q\in\mathbf{\s}} q\lcr \bigvee
R_\beta[q]\rcr\\[4mm]
\hfill=\bigvee_\beta R\, .\end{aligned}\]

Finally, let $S\subseteq\mathbf{\s}$ with $R\subseteq S$, then
obviously $R_\beta[p]\subseteq S_\beta[p]$ for all
$p\in\mathbf{\s}$, hence $\bigvee R_\beta[p]\subseteq \bigvee
S_\beta[p]$. As a consequence,
\[ \bigvee_\beta R= \bigcup_{p\in\mathbf{\s}}p[\bigvee
R_\beta[p]]\subseteq \bigcup_{p\in\mathbf{\s}}p[\bigvee
S_\beta[p]]=\bigvee_\beta S.\]

(2) By Axiom P3, $(\bigvee_\l R)_\beta[p]\in\l_\beta$ for all
$p\in\mathbf{\s}$ and all $\beta\in\Omega$; whence $\bigvee_\beta
\lp \bigvee_\l R\rp=\bigvee_\l R$. Thus, $R^n\subseteq(\bigvee_\l R)^n=\bigvee_\l R$
, since $R\subseteq\bigvee_\l R$.

(3) Suppose that $R\in\ovee_\alpha\l_\alpha$. Then, by Definition
\ref{DefinitionSeparatedProduct}, $R_\beta[p]\in\l_\beta$ for all
$p\in\mathbf{\s}$ and all $\beta\in\Omega$. Hence
$R_\beta[p]=\bigvee R_\beta[p]$ in $\l_\beta$. As a consequence,
$\bigvee_\beta R=R$.

Suppose now that $\bigvee_\beta R=R$ for some $\beta\in\Omega$.
Let $q\in \mathbf{\s}$. Then
\[\begin{aligned}R_\beta[q]=\lp \bigvee_\beta R\rp_\beta[q]=
\lp \bigcup_{p\in\mathbf{\s}}p[\bigvee R_\beta[p]]\rp_\beta[q]=\lp
q[\bigvee R_\beta[q]]\rp_\beta[q]\\
=\bigvee R_\beta[q]\, .\end{aligned}\]
As a consequence, if $\bigvee_\beta R=R$ for all $\beta\in\Omega$,
then $R_\beta[q]\in\l_\beta$ for all $\beta\in\l_\beta$ and all
$q\in \mathbf{\s}$, hence $R\in\ovee_\alpha\l_\alpha$.
\end{proof}

\begin{remark}
Note that $\bigvee_\alpha R$ is not necessarily in $\l$. Take for
instance $\l=\l_1\owedge\l_2$ and $R=\{p,q,r\}$ with $p_1\ne r_1$,
$p_2\ne r_2$, $q_1=p_1$, and $q_2=r_2$. Then $\bigvee_\owedge
R=p_1\t (p_2\vee r_2)\cup (p_1\vee r_1)\t r_2$. Now,
$\bigvee_2 R=p_1\t (p_2\vee r_2)\cup r\ne \bigvee_\owedge
R$. But, $\bigvee_1(\bigvee_2 R)=p_1\t (p_2\vee r_2)\cup
(p_1\vee r_1)\t r_2= \bigvee_\owedge R$.
\end{remark}

\begin{theorem}\label{TheoremTrain}
Let $\{\l_\alpha\}_{\alpha\in\Omega}$ be a family of simple
closure spaces on $\s_\alpha$ and
$\mathbf{\s}=\prod_\alpha\s_\alpha$. If there is at most one
$\beta\in\Omega$ such that $\l_{\beta}\ne 2^{\s_{\beta}}$, then
$\owedge_\alpha\l_\alpha=\ovee_\alpha\l_\alpha=\{
R\subseteq\mathbf{\s}\ \lvert\ R_{\beta}[p]\in\l_{\beta},\, \forall
p\in\mathbf{\s}\}$.\end{theorem}

\begin{proof} Let $R\subseteq\mathbf{\s}$.
Since $\owedge_\alpha\l_\alpha\subseteq\ovee_\alpha\l_\alpha$,
$\bigvee_\ovee R\subseteq\bigvee_\owedge R$. We prove that
$\bigvee_\owedge R\subseteq\bigvee_\ovee R$.
Write $\s_1=\s_{\beta}$, $\l_1=\l_{\beta}$, and
$\s_2=\prod_{\alpha\ne\beta}\s_\alpha$. Hence
$\mathbf{\s}=\s_1\t\s_2$. We denote by
$\pi_i\colon\mathbf{\s}\to\s_i$ ($i=1,2$) the $i-$th coordinate
map. From definition \ref{DefinitionSeparatedProduct}, we find
that $\owedge_{\alpha\ne\beta}\l_\alpha=2^{\s_2}$.
Therefore, since $\owedge$ is associative, $\mathop\owedge_{\alpha\in\Omega}\l_\alpha=\l_1\owedge 2^{\s_2}$.

Let $a\in\l_1$ and $b\in 2^{\s_2}$. 
Note that $R\subseteq a\t
\s_2\bigcup \s_1\t b$ if and only if $g(b):=\pi_1(\s_1\t b^c\cap R)\subseteq a$,
where $b^c$ denotes the set-complement of $b$, {\it i.e.},
$b^c=\s_2\backslash b$. 
Note also that 
\[R\subseteq 
\lp\vee g(b)\rp\t\s_2\cup\s_1\t \lp b\cap \pi_2(R)\rp\mbox{, for all }b\in2^{\s_2}.\]
Moreover, let $a\in\l_2$ and $b\in 2^{\s_2}$. Then, $R\subseteq a\t\s_2\cup\s_1\t b$ implies
\[\lp\vee g(b)\rp\t\s_2\cup\s_1\t \lp b\cap \pi_2(R)\rp
\subseteq a\t\s_2\cup\s_1\t b\ .\]
Finally, $g(b)=g(b\cap\pi_2(R))$ for all $b\in2^{\s_2}$.
As a consequence, we find that
\[\begin{aligned}\bigvee_\owedge
R&=\bigcap\lac a\t\s_2\cup\s_1\t b\ \lvert\ a\in\l_1,\,b\in2^{\s_2}\,\mbox{and}\,R\subseteq a\t\s_2\cup\s_1\t b\rac\\
&=\bigcap_{b\subseteq\pi_2(R)}\lp\bigvee g(b)\rp\t \s_2\cup \s_1\t b\\
&=\bigcup_{f\in 2^{\lp 2^{^{\pi_2(R)}}\rp}}\lp\vbox{\vspace{0.5cm}}\bigcap\lac\bigvee g(b)\ \lvert\ b\in
f^{-1}(1)\rac\rp\t\lp\bigcap f^{-1}(2)\rp\, .\end{aligned}\]
Write $X_f$ for the last term in the preceding equation and 
$b=\cap f^{-1}(2)$. Let $g\in2^{\lp 2^{^{\pi_2(R)}}\rp}$ such that
$b\subseteq c\,\Rightarrow\, g(c)=2$. Then
$\cap g^{-1}(2)=b$ and $X_f\subseteq X_g$. As a consequence, if for $b\subseteq \pi_2(R)$ we define
$m(b):=\{c\subseteq\pi_2(R)\ \lvert\ b\nsubseteq c\}$, then,
\[\bigvee_\owedge R=\bigcup_{b\subseteq\pi_2(R)}\lp\vbox{\vspace{0.5cm}}\bigcap\lac\bigvee
g(c)\ \lvert\ c\in m(b)\rac\rp\t b.\]

Note that for all $q\in b$, $q^c\bigcap\pi_2(R)\in m(b)$.
Moreover,
\[g\lp q^c\bigcap\pi_2(R)\rp=\pi_1\lp\vbox{\vspace{0.5cm}}\s_1\t\lp q\bigcup \pi_2(R)^c\rp\bigcap
R\rp=R_1[(\cdot,q)].\]
As a consequence,
\[\begin{aligned}\bigvee_\owedge R&\subseteq
\bigcup_{b\subseteq\pi_2(R)}\lp\vbox{\vspace{0.5cm}}\bigcap\lac\bigvee
R_1[(\cdot,q)]\ \lvert\ q\in b\rac\rp\t
b\\
&=\bigcup_{q\in \pi_2(R)}\lp\bigvee R_1[(\cdot,q)]\rp\t
q=\bigvee_{\beta} R\, .\end{aligned}\]
Finally, by Lemma \ref{LemmaJoininTensorF}, $\bigvee_{\beta}
R=\bigvee_\ovee R$. As a consequence, $\bigvee_\owedge
R\subseteq\bigvee_\ovee R$.
\end{proof}

\begin{definition}\label{DefinitionXi}
Let $\s$ be a nonempty set. We denote by $\mo_\s$ the simple
closure space on $\s$ which contains only $\emptyset$, $\s$, and
all singletons of $\s$. We write $\mo_n$ if $\abs\s=n$.

Let $\l$ be a simple closure space. We say that $\l$ {\it contains $\mo_n$} if there
are $n$ atoms $p_1,\cdots,p_n$ such that $p_1\bigvee p_n\gtrdot
p_i$ for all $i$ between $1$ and $n$.

Let $\{\s_\alpha\}_{\alpha\in\Omega}$ be a family of sets, and
$\mathbf{\s}=\prod_\alpha\s_\alpha$. We denote by
$\Xi(\mathbf{\s})$ the set
$\{R\subseteq\mathbf{\s}\ \lvert\ p_\alpha\ne q_\alpha,\,\forall p,q\in
R,\,\alpha\in\Omega\}$.\end{definition}

\begin{theorem}\label{TheoremPerspectiveThenAertsdiffOtf}
Let $\l_1$ and $\l_2$ be simple closure spaces on $\s_1$ and
$\s_2$ respectively. If both $\l_1$ and $\l_2$ contain $\mo_3$,
then $\l_1\owedge\l_2\ne\l_1\ovee\l_2$.
\end{theorem}

\begin{proof}
Let $R=\{p,q,r\}\in\Xi(\s_1\t\s_2)$ such that for $i=1$ and for
$i=2$, $p_i\vee q_i$ covers $p_i$, $q_i$ and $r_i$. By
Definition \ref{DefinitionSeparatedProduct},
$R=\{p,q,r\}\in\l_1\ovee\l_2$. On the other hand, $R\subseteq a\t
\s_2\cup \s_1\t b$ if and only if $p_1\vee q_1\subseteq a$
or $p_2\vee q_2\subseteq b$. As a consequence, 
$\bigvee_\owedge R=\lp p_1\vee q_1\rp\t \lp p_2\vee q_2\rp\ne R=\bigvee_\ovee
R$.
\end{proof}

The next corollary is a partial converse to Theorem
\ref{TheoremTrain}.

\begin{corollary}
Let $\{\l_i\}_{1\leq i\leq n}$ be a finite family of simple
closure spaces on $\s_i$. Suppose that each $\l_i$ is a
DAC-lattice. If $\owedge_{i=1}^n\l_i=\ovee_{i=1}^n\l_i$, then there
is at most one $i$ such that $\l_i\ne 2^{\s_i}$.
\end{corollary}

\begin{proof} Suppose that $\l_i\ne 2^{\s_i}$ for $i=m$ and for $i=k$
with $m\ne k$ between $1$ and $n$, and that
$\owedge_{i=1}^n\l_i=\ovee_{i=1}^n\l_i$. Then, for $i=m$ and for
$i=k$, there is an atom $p_i$ which is not a central element.
Therefore, there is an atom $q_i$ such that $p_i\bigvee q_i$
contains a third atom (see \cite{Maeda/Maeda:handbook}, Theorems
28.8, 27.6 and Lemma 11.6).

Let $\l_0:=\owedge_{i\ne m}\l_i$ and $\l_1:=\ovee_{i\ne m}\l_i$. By
hypothesis, $\l_m\owedge\l_0=\l_m\ovee\l_1$. Therefore, from Lemma
\ref{LemmaAlex} and Axiom P3, $\l_0=\l_1$. Let $r\in\prod_{i\ne
m}\s_i$ be an atom of $\l_0$. From Lemma \ref{LemmaAlex},
$r[p_k]\vee_\owedge r[q_k]=r[p_k\vee q_k]$, hence contains
a third atom, therefore $\l_0$ contains $\mo_3$. As a consequence,
from Theorem \ref{TheoremPerspectiveThenAertsdiffOtf}, we find
that
\[\l_m\ovee\l_1=\l_m\owedge\l_0\ne\l_m\ovee\l_0=\l_m\ovee\l_1,\]
a contradiction.
\end{proof}

\begin{theorem}\label{TheoremCoatomsCountable}
Let $\l$ be a coatomistic simple closure space on $\s$. Suppose
that for any countable set $A$ of coatoms $\l$, $\bigcup A\ne\s$.
For all integers $i$, let $\l_i=\l$. Then
$\owedge_{i=1}^\infty\l_i\ne\ovee_{i=1}^\infty\l_i$.
\end{theorem}

\begin{proof} For all $i$, let $\s_i=\s$.
Denote $\prod_{i=1}^\infty\s_i$ by $\mathbf{\s}$. Let
\[R=\{p\in\mathbf{\s}\ \lvert\ p_n=p_m,\,\forall m,n\}\ne\mathbf{\s}.\]
From Definition \ref{DefinitionSeparatedProduct},
$R\in\ovee_{i=1}^\infty\l_i$. On the other hand, by hypothesis, for
any $a\in\prod_{i=1}^\infty(\l_i\backslash\{1\})$, $R\nsubseteq\bigcup_{i=1}^\infty
\pr_i(a_i)$, hence $\bigvee_\owedge R=1$.\end{proof}

\begin{example} The Hypothesis of Theorem
\ref{TheoremCoatomsCountable} is fulfilled for instance if
$\l=\pro(\h)$ with $\h$ a real or complex Hilbert
space.\end{example}
\section{Orthocomplementation}\label{SectionOrthocomplementation}

This section is devoted to our main results, which show, subject
to weak conditions, that if $\l\in\sep(\l_1,\cdots,\l_n)$ is
orthocomplemented, then $\l=\owedge_{i=1}^n\l_i$. Using Theorem
\ref{TheoremAutoFactorwithWeaklyConnected}, we prove that this
holds if each $\l_i$ is weakly connected (hence, irreducible) and
coatomistic and if $\l$ is moreover transitive. For the second result, we assume that
each $\l_i$ is orthocomplemented and that all its irreducible
components different from $2$ are weakly connected, but we do not
need to assume that $\l$ is transitive.

\begin{definition} Let $\{\s_\alpha\}_{\alpha\in\Omega}$ be a
family of nonempty sets, $p\in\mathbf{\s}=\prod_\alpha\s_\alpha$,
$\beta,\,\gamma\in\Omega$, $A_\beta\subseteq\s_\beta$, and
$A_\gamma\subseteq\s_\gamma$. Then
\[ p[A_\beta,A_\gamma]:=\{q\in\mathbf{\s}\ \lvert\ 
q_\beta\in A_\beta,\,q_\gamma\in A_\gamma,\,
q_\alpha=p_\alpha,\,\forall \alpha\ne \beta,\, \gamma\}.\]
\end{definition}

\begin{lemma}\label{LemmaCoatomsAertsCoatomsSep}
Let $\{\l_i\}_{1\leq i\leq n}$ be a finite family of simple
closure spaces on $\s_i$, $x\in\prod_{i=1}^n\l_i$ with $x_i$
coatoms, and let $\l\in\sep(\l_i,1\leq i\leq n)$. Then:
\begin{enumerate}
\item $X:=\bigcup_{i=1}^n \pr_i(x_i)$ is a coatom of $\l$.
\item Let $a_j\in\l_j$ for some $j$ between 1 and $n$ and let $Z$ be a coatom
of $\l$ above $\bigcup_{i\ne j} \pr_i(x_i)\bigcup\pr_j(a_j)$. Then there is a coatom $z_j$ 
of $\l_j$ such that $Z=\bigcup_{i\ne j} \pr_i(x_i)\bigcup\pr_j(z_j)$.
\end{enumerate}
\end{lemma}

\begin{proof}
(1): By Axiom P2, $X\in\l$. Let $\mathbf{\s}=\prod_{i=1}^n\s_i$ and
$p\in\mathbf{\s}$ not in $X$. Write $R^0:=p\cup X$ and
$y:=p\vee X$. Define $R^N:=\bigvee_n\cdots\bigvee_1 R^{N-1}$.
By Lemma \ref{LemmaJoininTensorF}, $R^N\subseteq y$, for all $N$.
Now, 
\[
R_k^0[q]=\lac s\in\s_k\ \lvert\ q[s,k]\in p\cup X\rac=\{p\}_k[q]\cup  X_k[q];
\]
therefore, 
\[
q[R_k^0[q]]=q[X_k[q]]\mbox{ or }q[R_k^0[q]]=p[x_k\cup p_k] ,
\]
hence 
\[
q[\bigvee R_k^0[q]]=q[\vee X_k[q]]\subseteq X\mbox{ or }q[\vee R_k^0[q]]=p[\s_k]\,.
\]
As a consequence, 
$\bigvee_k R^0= p[\s_k]\cup R^0$, thus
$R^1=R^0\bigcup_{k=1}^n p[\s_k]$. Further,
\[R^2=R^0\bigcup\lac p[\s_{k_1},\s_{k_2}]\ \lvert\ 1\leq k_1\leq n-1\mbox{ and }k_1+1\leq k_2\leq n\rac .\]
Hence obviously, for $N=n$, $R^N=\mathbf{\s}$.

(2): Write $Y:=\bigcup_{i\ne j} \pr_i(x_i)\bigcup\pr_j(a_j)$. 
Let $r\in Z$ not in $Y$. Note that
$r[a_j]\subseteq Y$. Therefore, by Lemma
\ref{LemmaAlex}, $r[a_j\bigvee r_j]\subseteq Z$.
Let $k\ne j$ between $1$ and $n$. Then,
$r[x_k,a_j\vee r_j]\subseteq Y$. Hence,
since $r[a_j\vee r_j]\subseteq Z$, by Lemma
\ref{LemmaAlex}, $r[\s_k,a_j\vee r_j]\subseteq
Z$. Repeating this argument, we find that
$\pr_j(a_j\vee r_j)\subseteq Z$. As a consequence,
\[
Z=Y\bigvee_{r\in Z\backslash Y}r
\subseteq \bigcup_{i\ne j}\pr_i(x_i)\bigvee_{r\in Z\backslash Y}\pr_j(a_j\vee r_j)\subseteq Z
\]
which proves that $Z=\bigcup_{i\ne j}\pr_i(x_i)\cup \pr_j(z_j)$ for some coatom $z_j$ of $\l_j$.
\end{proof}

\begin{proposition} Let $\{\l_\alpha\}_{\alpha\in\Omega}$ be a
family of simple closure spaces on $\s_\alpha$.
\begin{enumerate}
\item If $\owedge_\alpha\l_\alpha$ is coatomistic, then all
$\l_\alpha$'s are coatomistic.
\item Suppose that $\Omega$ is finite and that all $\l_\alpha$'s
are weakly connected. If $\owedge_\alpha\l_\alpha$ is
orthocomplemented, then all $\l_\alpha$'s are orthocomplemented.
\end{enumerate}
\end{proposition}

\begin{proof}
(1): Let $\mathbf{\s}=\prod_\alpha\s_\alpha$ and $p\in\mathbf{\s}$.
From Definition \ref{DefinitionSeparatedProduct}, a coatom of
$\owedge_\alpha\l_\alpha$ over $p$ is necessarily of the form
$\bigcup_\alpha\pr_\alpha(x_\alpha)$ where all $x_\alpha$'s are
coatoms and with $p_\beta\in x_\beta$ for some $\beta\in\Omega$.
Let $\s'_\alpha[p_\alpha]$ denote the set of coatoms of $\l_\alpha$ above
$p_\alpha$ (write $\s_\alpha':= \s_\alpha'[0_\alpha]$).
Define $\chi=\lac x\in\prod_\alpha\s'_\alpha\ \lvert\ p_\beta\in x_\beta\
\mbox{for\ some}\ \beta\in\Omega\rac$.
If
$\owedge_\alpha\l_\alpha$ is coatomistic, then
\[p=\bigcap\lac\bigcup_\alpha\pr_\alpha(x_\alpha)\ \lvert\ 
x\in\ \chi\rac=\bigcup_{f\in\Omega^{^\chi}}\lp\vbox{\vspace{0.5cm}}
\prod_{\alpha\in\Omega}\lp\bigcap f^{-1}(\alpha)_\alpha\rp\rp.\]
Hence, if $\s'_{\alpha_0}[p_{\alpha_0}]=\emptyset$ for some $\alpha_0\in\Omega$, 
then $\pi_{\alpha_0}(p)=
\s_{\alpha_0}$ which means that $\l_{\alpha_0}=2$.
As a consequence, all $\l_\alpha$'s are coatomic
and we can assume that $\s'_\alpha[p_\alpha]\ne\emptyset$ for all
$\alpha\in\Omega$. Moreover, it follows that 
$p_\alpha=\bigcap\lac x_\alpha\in\s'_\alpha[p_\alpha]\rac$,
for all $\alpha\in\Omega$, hence that all
$\l_\alpha$'s are coatomistic.

(2): Let $\beta\in\Omega$ and $a\in\l_\beta$. We first prove that
there is $b\in\l_\beta$ such that $\pr_\beta(a)'=\pr_\beta(b)$.
Let $f\colon\mathbf{\s}\to\prod\s_\alpha'$ such that
$p'=\bigcup_\alpha\pr_\alpha(f(p)_\alpha)$ for all
$p\in\mathbf{\s}$. Note that $f$ is injective.

(2.1) {\bf Claim}: Let $p,\, q\in\mathbf{\s}$.  If $f(p)_\alpha\ne
f(q)_\alpha$ for at least two $\alpha\in\Omega$, then $p\vee
q=p\cup q$. 

{\it Proof}\nobreak : Let
$\Omega_{\ne}=\{\alpha\in\Omega\ \lvert\ f(p)_\alpha\ne f(q)_\alpha\}$
and $\gamma\ne\delta\in\Omega_{\ne}$. Let $r\in p\vee q$. Then,
\[\begin{aligned}
\bigcup_\alpha\pr_\alpha(f(r)_\alpha)=r'\supseteq p'\bigcap q'\supseteq\hspace{7cm}\\
\lp\pr_\gamma(f(p)_\gamma)\cap\pr_\delta(f(q)_\delta)\rp\cup
\lp\pr_\delta(f(p)_\delta)\cap\pr_\gamma(f(q)_\gamma)\rp \, .
\end{aligned}\]
Using the notation of Remark \ref{remarkread}, we can write for $\Omega=\{1,2,3\}$,
$\gamma=1$ and $\delta=2$,
\[f(r)_1\vert\vert\cup\vert f(r)_2\vert\cup\vert\vert f(r)_3\supseteq 
f(p)_1f(q)_2\vert\cup f(q)_1f(p)_2\vert\]
Therefore, $f(r)_\gamma=x_\gamma$ and $f(r)_\delta=x_\delta$ with
$x=f(p)$ or $x=f(q)$. As a consequence, $r'=p'$ or $r'=q'$, hence
$r=p$ or $r=q$, proving the claim.
\begin{figure}[t]
\setlength{\unitlength}{1cm} 
\begin{picture}(10,3.3)(0,0)
\put(9.5,1.5){\makebox(0,0){$p$}}
\put(1,1.5){\linethickness{0.05 cm}\line(1,0){1.9}}
\put(3.1,1.5){\linethickness{0.05 cm}\line(1,0){5.9}}
\put(3,0.05){\line(0,1){0.3}}
\put(3,0.5){\line(0,1){0.3}}
\put(3,1){\line(0,1){0.3}}
\put(3,1.5){\line(0,1){0.3}}
\put(3,2){\line(0,1){0.3}}
\put(3,2.5){\line(0,1){0.3}}
\put(3,0.5){\circle*{0.2}}
\put(3.5,0.5){$q_\gamma$}
\put(3,2.5){\circle*{0.2}}
\put(3.5,2.5){$p_\gamma$}
\put(2.9,3){$\gamma$}
\put(2.8,1.7){\linethickness{0.03 cm}\line(0,1){0.9}}
\put(2.3,2){$A_1$}
\put(2.8,0.4){\linethickness{0.03 cm}\line(0,1){0.9}}
\put(2.3,0.8){$A_3$}
\put(3.2,1){\linethickness{0.03 cm}\line(0,1){1}}
\put(3.3,1){$A_2$}
\multiput(2.8,2)(0.05,0){8}{\line(1,0){0.025}}
\multiput(2.8,1.7)(0.05,0){8}{\line(1,0){0.025}}
\multiput(2.8,1.3)(0.05,0){8}{\line(1,0){0.025}}
\multiput(2.8,1)(0.05,0){8}{\line(1,0){0.025}}
\end{picture}
\caption{} \label{fig1}
\end{figure}

(2.2) {\bf Claim}: Let $p,\, q\in\mathbf{\s}$ such that $p$ and
$q$ differ only by one component, say $\gamma$. Then
$f(p)_\alpha=f(q)_\alpha$ for all $\alpha\ne\gamma$. 

{\it
Proof}\nobreak : Since $\l_\gamma$ is weakly connected, from part
2.1 we find that there is $\delta\in\Omega$ such that for all
$r\in\s_\gamma$; we have $f(p[r,\gamma])_\alpha=f(p)_\alpha$ for all
$\alpha\ne\delta$ (see Figure \ref{fig1}).
Hence,
\[\begin{aligned}
p[\s_\gamma]'&=\bigcap_{r\in\s_\gamma}\lp\bigcup_{\alpha\ne\delta}
\pr_\alpha(f(p)_\alpha)\cup\pr_\delta(f(p[r,\gamma]_\delta)\rp\\
&=\bigcup_{\alpha\ne
\delta}\pr_\alpha(f(p)_\alpha)\cup\pr_\delta(x)\, ,
\end{aligned}\]
for some $x\in\l_\delta$. Now $p[\s_\gamma]\bigcap
\pr_\gamma(f(p)_\gamma)\ne 0$. Therefore, if $\delta\ne\gamma$,
then $p[\s_\gamma]\bigcap p[\s_\gamma]'\ne 0$, a contradiction. As
a consequence, $\delta=\gamma$ and $x=0$, proving the claim.

Let $p,\,q\in\mathbf{\s}$ such that $p\in\pr_\beta(a)'$, for some $a\in\l_\beta$, and such
that $p$ and $q$ differ only by one component, say $\gamma$,
with $\gamma\ne\beta$. Then $\pr_\beta(a)\subseteq
p'=\bigcup_\alpha\pr_\alpha(f(p)_\alpha)$, hence $a\subseteq
f(p)_\beta$. Now, from part 2.2, $f(q)_\beta=f(p)_\beta$,
therefore $\pr_\beta(a)\subseteq q'$, thus $q\in\pr_\beta(a)'$.
As a consequence, for all
$q\in\mathbf{\s}$ such that $q_\beta=p_\beta$,
$q\in\pr_\beta(a)'$.

Thus, we have proved that there is an element in
$\l_\beta$, which we denote by $a^{\p_\beta}$, such that
$\pr_\beta(a)'=\pr_\beta(a^{\p_\beta})$. Obviously, the mapping
$\p_\beta\colon\l_\beta\to\l_\beta$ is an orthocomplementation.
\end{proof}

\begin{theorem}
Let $\{\l_i\}_{1\leq i\leq n}$ be a finite family of coatomistic
weakly connected simple closure spaces on $\s_i$, and let
$\l\in\sep(\l_i,1\leq i\leq n)$. If $\l$ is transitive and
orthocomplemented, then $\l=\owedge_{i=1}^n\l_i$.
\end{theorem}

\begin{proof} Let $\mathbf{\s}=\prod_{i=1}^n\s_i$. We denote the
orthocomplementation of $\l$ by $'$. Let $x\in\prod_{i=1}^n\l_i$
with $x_i$ coatoms, and $X:=\bigcup_{i=1}^n \pr_i(x_i)$. By Lemma
\ref{LemmaCoatomsAertsCoatomsSep}, $X'=p$, for some
$p\in\mathbf{\s}$. Let $q\in\mathbf{\s}$. Since $\l$ is
transitive, there is $u\in\aut(\l)$ such that $u(p)=q$, hence
$q'=u(p)'$. Define $u'\in\aut(\l)$ as $u'(a):=u(a')'$. Then
$q'=u(p)'=u'(p')=u'(X)$. From Theorem
\ref{TheoremAutoFactorwithWeaklyConnected}, $u'$ factors,
therefore there is $y\in\prod_{i=1}^n\l_i$ with $y_i$ coatoms such
that $q'=Y:=\bigcup_{i=1}^n \pr_i(y_i)$.
\end{proof}

\begin{remark}
If
$\l\in\Sep_T(\l_i,1\leq i\leq n)$, $\l_i$ are transitive, and the
$u$ in Axiom P4 is an ortho-isomorphism of $\l$ for all $v_i\in
T_i$, then the proof does not require Theorem
\ref{TheoremAutoFactorwithWeaklyConnected}. Below we give a second
proof which requires neither Theorem
\ref{TheoremAutoFactorwithWeaklyConnected} nor that $\l$ be
transitive.
\end{remark}

\begin{definition}\label{defvac}
Let $\{\l_\alpha\}_{\alpha\in\Omega}$ be a family of
orthocomplemented simple closure spaces on $\s_i$. For 
$p\in\prod_{\alpha\in\Omega}\s_\alpha$ define $p^\#=\bigcup_{\alpha\in\Omega}\pr_\alpha(p_\alpha^{\p_\alpha})$
\end{definition}

\begin{remark}
It is easy to check that the mapping $p\mapsto p^\#$ is an orthocomplementation of 
$\owedge_{\alpha\in\Omega}\l_\alpha$ (see \cite{Ischi:2007} for details).
\end{remark}

\begin{remark}
Let $\l$ be an orthocomplemented simple closure space on $\s$.
Let $\mathcal{Z}(\l)$ denotes the center of $\l$. Then
$\mathcal{Z}(\l)=\{a\in\l\ \lvert\ a^\p=a^c\}$, where $a^c:=\s\backslash a$. 
For $p\in\s$, we write $e(p)$ for the central cover of $p$, that is 
$e(p)=\bigcap\{a\in\mathcal{Z}(\l)\ \lvert\ p\in a\}$.
\end{remark}

\begin{theorem}
Let $\{\l_i\}_{1\leq i\leq n}$ be a finite family of
orthocomplemented simple closure spaces on $\s_i$. Suppose
moreover that for each $i$, all irreducible components of $\l_i$
different from $2$ are weakly connected. Let $\l\in\sep(\l_i,1\leq
i\leq n)$. If $\l$ is orthocomplemented, then
$\l=\owedge_{i=1}^n\l_i$.
\end{theorem}

\begin{proof} Let $\mathbf{\s}=\prod_{i=1}^n\s_i$. We denote
the orthocomplementation of $\l$ by $'$ and the orthocomplementation of $\l_i$ by $\p_i$. 
Since coatoms
of $\owedge_i\l_i$ are coatoms of $\l$, we can define a map $\phi$ on
$\mathbf{\s}$ as $\phi(p)=p^{\#'}$. Note that $\phi$ is injective. We
prove in four steps that $\phi$ is surjective.

(1) {\bf Claim}: Let $p\in\mathbf{\s}$ and $a_j\in\l_j$ for some $j$ between
$1$ and $n$; then  $\phi(p[a_j])=p[a_j]^{\#'}$. 

{\it
Proof}\nobreak : If $q\in p[a_j]$, then $p[a_j]^\#\subseteq q^\#$,
hence $\phi(q)\in p[a_j]^{\#'}$. On the other hand, if $q\in
p[a_j]^{\#'}$, then, $p[a_j]^\#=\bigcup_{i\ne j} \pr_i(p_i^{\p_i})\cup
\pr_j(a_j^{\p_j})\subseteq
q'$.

Thus, by Lemma \ref{LemmaCoatomsAertsCoatomsSep}, there is an atom $s_j$ with $a_j^{\p_j}\subseteq s_j^{\p_j}$ 
and 
\[q'=\bigcup_{i\ne j}\pr_i(p_i^{\p_i})\bigcup \pr_j(s_j^{\p_j}).\]
Therefore, $q=p[s_j]^{\#'}$ with $s_j\in a_j$. Hence
$q\in\phi(p[a_j])$, proving the claim.

(2) {\bf Claim}: Let $p\in\mathbf{\s}$.  For all $j$ between $1$
and $n$ there is $q\in \prod_{i=1}^n e(p_i)$ and $k$ between $1$
and $n$ such that $\phi(p[e(p_j)])\subseteq q[e(p_k)]$. 

{\it Proof}\nobreak
: Note that if $a\in\l$, then $a'\subseteq a^c$, where $a^c$
denotes the set-complement of $a$, {\it i.e.},
$a^c=\mathbf{\s}\backslash a$. By Claim 1,
\[\phi(p[e(p_j)])\subseteq (p[e(p_j)]^\#)^c=e(p_j)
\mbox{$\prod$}_{i\ne j}
((p_i^{\p_i})^c)\subseteq\mbox{$\prod$}_{i=1}^n e(p_i),\]
since if $q_i\in e(p_i)$, then $(q_i^{\p_i})^c\subseteq e(p_i)$.
If $[0,e(p_j)]=2$ ({\it i.e.}, $e(p_j)=p_j$), then the proof of Claim 2 is
finished.

Otherwise, let $t_j\ne s_j\in e(p_j)$. Then,
\[\begin{aligned}q\in\phi(p[t_j])\vee\phi(p[s_j])&\mbox{ iff }
p[t_j]^\#\cap p[s_j]^\#\subseteq q'\\
&\mbox{ iff } \bigcup_{i\ne j}\pr_i(p_i^{\p_i})\cup
\pr_j((t_j\bigvee s_j)^{\p_j}) \subseteq q'\, ,\end{aligned}\]
hence, by Lemma \ref{LemmaCoatomsAertsCoatomsSep}, if and only if 
$q=\phi(p[r_j])$ for some $r_j\in t_j\vee
s_j$.

As a consequence, if $t_j\bigvee s_j$ contains a third atom, so
does $\phi(p[t_j])\vee \phi(p[s_j])$. Hence, by Lemma
\ref{LemmaJoinP1diffQ1andP2diffQ2} part 1, $\phi(p[t_j])$ and
$\phi(p[s_j])$ differ only by one component, say $k$. Therefore, 
since $\l_j$ is weakly connected, for all
$A_j^\gamma$ in the connected covering of $\s_j$, there is
$q^\gamma\in \prod_{i=1}^n e(p_i)$ and $k_\gamma$ such that
$\phi(p[e(p_j)\bigcap A_j^\gamma])\subseteq
q^\gamma[e(p_{k_\gamma})]$. From Hypotheses 1 and 3 in Definition
\ref{DefinitionConnected}, the maps $\gamma\mapsto k_\gamma$ and
$\gamma\mapsto q_\gamma$ are constant since $\phi$ is injective. This completes
the proof of the claim.

(3) {\bf Claim}: For all $p\in\mathbf{\s}$ and all $j$ between $1$
and $n$, there is $q\in\prod_{i=1}^n e(p_i)$ such that
$\phi(p[e(p_j)])=q[e(p_j)]$. 

{\it Proof}\nobreak : From Claim 2,
there is $q\in\prod_{i=1}^n e(p_i)$, $k$ between $1$ and $n$, and
$b_k\subseteq e(p_k)$, such that $\phi(p[e(p_j)])= q[b_k]$.

Assume first that $k\ne j$. Let $R^0:=q[b_k]\cup p[e(p_j)]^\#$.
By part 1,
\[{R^0}'=q[b_k]'\cap p[e(p_j)]^{\#'}=\phi(p[e(p_j)])'\cap
\phi(p[e(p_j)])=0,\]
hence $\bigvee R^0=1$. On the other hand, note that since $q_j\in e(p_j)$,
$q_j\vee e(p_j)^c=q_j\cup e(p_j)^c$. If $e(p_j)=p_j$, then part 3 is trivial. 
Hence we can assume that $e(p_j)\ne p_j$, thus that $q_j\vee e(p_j)^c\ne 1$, in other words that 
$e(p_j)\cap q_j^{\p_j}\ne\emptyset$. Now, for any $r_j\in e(p_j)\cap q_j^{\p_j}$, we have
\[
R^0=q[b_k]\bigcup_{i\ne j} \pr_i(p_i^{\p_i})\cup
\pr_j(e(p_j)^c)\subseteq \pr_j(r_j^{\p_j})\bigcup_{i\ne
j}\pr_i(p_i^{\p_i})\, ,
\]
whence by Axiom P2,
$\bigvee R^0\ne 1$, a contradiction. As a consequence, $k=j$.

Let $R:=q[b_j]\bigcup p[e(p_j)]^\#$. From Claim 1,
$R'=\phi(p[e(p_j)])'\cap \phi(p[e(p_j)])=0$. Therefore,
$\bigvee R=1$. Now, $R\subseteq\bigcup_{i\ne j}\pr_i(p_i^{\p_i})\cup
\pr_j(e(p_j)^c\vee b_j)$.
Note that for all $a\in\l_i$ with $a\subseteq e(p_j)$, we have
$e(p_j)\cap(a\vee e(p_j)^{\p_j})=a$. Therefore, $a\vee
e(p_j)^{\p_j}=a\cup e(p_j)^{\p_j}$. Hence, we find that
\[R\subseteq\bigcup_{i\ne j}\pr_i(p_i^{\p_i})\cup
\pr_j(e(p_j)^c\cup b_j).\]
Since $\bigvee R=1$, from Axiom P2 we find that
$b_j=e(p_j)$, proving the claim.

(4): Let $p\in\mathbf{\s}$ and $s\in \prod_{i=1}^n e(p_i)$. By Claim
3, $\phi(p[e(p_1)])=q^1[e(p_1)]$. Therefore, there is $r_1\in
e(p_1)$ such that $\phi(p[r_1])_1=s_1$. Let $k\leq n$ and $r_1\in
e(p_1), \cdots,\,r_k\in e(p_k)$ such that
$\phi(p[r_1,\cdots,r_k])_i=s_i$, for all $i\leq k$, and such that
$\phi(p[r_1,\cdots,r_k])_{k+1}$ is different from $s_{k+1}$. By
Claim 3,
\[\phi(p[r_1,\cdots,r_k,e(p_{k+1})])=q^{k+1}[s_1,\cdots
s_k,e(p_{k+1})].\]
Hence there is $r_{k+1}\in e(p_{k+1})$ such that
$\phi(p[r_1,\cdots,r_{k+1}])_i=s_i$, for all $i$ between $1$ and
$k+1$. As a consequence, $\phi$ is
surjective.
\end{proof}
\section{Covering property}\label{SectionCoevringProperty}

In this Section, we prove, under some assumptions, that the top element
$\ovee_\alpha\l_\alpha$ has the covering property if and only at
most one $\l_\alpha$ is not a power set. We reproduce the
analogue result concerning the bottom element
$\owedge_\alpha\l_\alpha$ which is due to Aerts \cite{Aerts:1982}.
Moreover, for $\l_i=\mo_{\s_i}$ ($i=1,2$) and
$T=\aut(\l_1)\t\aut(\l_2)$, we prove that there is a unique
$\l\in\Sep_T(\l_1,\l_2)$ with the covering property.

\begin{theorem}[D. Aerts, \cite{Aerts:1982}]
Let $\{\l_\alpha\}_{\alpha\in\Omega}$ be a family of orthocomplemented simple
closure spaces on $\s_\alpha$. If $\owedge_\alpha\l_\alpha$ has the
covering property or is orthomodular, then there is at most one
$\beta\in\Omega$ such that $\l_{\beta}\ne 2^{\s_{\beta}}$.
\end{theorem}

\begin{proof} Let $\l$ be an orthocomplemented 
simple closure space on $\s$ and let
$p,\, q\in\s$ such that $p\vee q=p\cup q$. Define
$x:=q^\p\cap(p\vee q)$, then $x=0$ or $x=p$. If
$\l$ has the covering property, then $1\ne x^\p$, whereas if $\l$
is orthomodular, $x\vee q=p\vee q$. As a consequence, $x=p$,
hence $p\p q$.

Let $\beta\in\Omega$. Suppose that $\l_{\beta}\ne 2^{\s_{\beta}}$.
Then there are two non orthogonal atoms, say $r_{\beta}$ and
$s_{\beta}$. Let $r_\gamma,\,s_\gamma\in\s_\gamma$ for some
$\gamma\in\Omega$ different from $\beta$, and let $p,\,
q\in\prod_\alpha\s_\alpha$ defined as $p_\alpha=q_\alpha$, for all
$\alpha\ne \beta,\,\gamma$, and $p_\alpha=r_\alpha$ and $q_\alpha=
s_\alpha$ if $\alpha=\beta$ or $\gamma$. By Lemma
\ref{LemmaJoinP1diffQ1andP2diffQ2}, $p\vee q=p\cup q$.
Therefore, by what precedes, since $\#$ is an orthocomplementation
of $\owedge_\alpha\l_\alpha$, so $p\# q$, hence by Definition \ref{defvac}, 
$r_\gamma\p_\gamma s_\gamma$. As a
consequence, $\l_\gamma=2^{\s_\gamma}$.\end{proof}

\begin{proposition} Let $\{\l_\alpha\}_{\alpha\in\Omega}$ be a
family of simple closure spaces on $\s_\alpha$.
If $\ovee_\alpha\l_\alpha$ has the covering property, then all
$\l_\alpha$'s have the covering property.\end{proposition}

\begin{proof} Let $a_\beta\in\l_\beta$, $q_\beta\in\s_\beta$
not in $a_\beta$, and $p\in\prod_\alpha\s_\alpha$. By the covering
property, we find that $p[q_\beta]\bigvee p[a_\beta]\gtrdot
p[a_\beta]$; Whence by Lemma \ref{LemmaAlex}.3, $q_\beta\bigvee
a_\beta\gtrdot a_\beta$.\end{proof}

\begin{remark}In the next theorem, we assume that each $\alpha$,
$\l_\alpha$ is a simple closure space on $\s_\alpha$, that $\l_\alpha\ne 2^{\s_\alpha}$,
and that $\l_\alpha$
contains $\mo_4$. This is for instance the case if $\l_\alpha$ is
orthocomplemented orthomodular with the covering property. Indeed,
if $\l_\alpha\ne 2^{\s_\alpha}$, there is an atom $p_\alpha$ which
is not central, hence such that $e(p_\alpha)$ contains at least
two atoms, say $r$ and $s$. Moreover, $r\vee s$ contains at
least three atoms, for $\l_\alpha$ has the covering property.
Finally, since $\l_\alpha$ is orthomodular, $[0,r\vee s]$ is
orthocomplemented, hence contains at least four atoms.
\end{remark}

\begin{theorem}\label{TheoremCoveringinTensorF}
Let $\{\l_\alpha\}_{\alpha\in\Omega}$ be a family of simple
closure spaces on $\s_\alpha$.
\begin{enumerate}
\item If each $\l_\alpha$ has the covering property and if there is
at most one $\beta\in\Omega$ such that $\l_\beta\ne 2^{\s_\beta}$,
then $\ovee_\alpha\l_\alpha$ has the covering property.
\item Suppose that each $\l_\alpha$ different from $2^{\s_\alpha}$
contains $\mo_4$. If $\ovee_\alpha\l_\alpha$ has the covering
property, then there is at most one $\beta\in\Omega$ such that
$\l_\beta\ne 2^{\s_\beta}$.
\end{enumerate}
\end{theorem}

\begin{proof} (1): Let $\mathbf{\s}=\prod_\alpha\s_\alpha$,
$a\in\ovee_\alpha\l_\alpha$, $q\in\mathbf{\s}$ not in $a$, and
$R=q\cup a$. By Lemma \ref{LemmaJoininTensorF}, $q\vee_\ovee a=\bigvee_{\beta}R=\bigcup_{p\in\mathbf{\s}}p[\bigvee
R_{\beta}[p]]$.
Now, by definition,
\[R_{\beta}[p]=\pi_\beta(p[\s_\beta]\cap(q\cup a))
=q_\beta[p]\cup a_\beta[p],\]
and $q_{\beta}[p]=q_{\beta}$ if $p\in q[\s_{\beta}]$, and
$q_{\beta}[p]=\emptyset$ otherwise. Hence, for all $p$ not in
$q[\s_{\beta}]$, $R_{\beta}[p]=a_{\beta}[p]\in\l_{\beta}$, and if
$p\in q[\s_{\beta}]$, then $R_{\beta}[p]=q_{\beta}\cup
a_{\beta}[p]$. Therefore, we find that
\[q\cup a\subseteq q\vee_\ovee a=
q[q_{\beta}\vee a_{\beta}[q]]\bigcup_{p\not\in
q[\s_{\beta}]}p[a_{\beta}[p]] \subseteq a\cup
q[q_{\beta}\bigvee a_{\beta}[q]].\]
As a consequence, $q\bigvee_{\ovee} a\gtrdot a$.

(2): Let $\beta\ne\gamma\in\Omega$. Suppose that neither
$\l_\beta=2^{\s_\beta}$ nor $\l_\gamma=2^{\s_\gamma}$. Let
$p,\,q,\,r,\,s,\,t\in\mathbf{\s}$ such that
$p_\alpha=q_\alpha=r_\alpha=s_\alpha=t_\alpha$, for all $\alpha$
different from $\beta$ and $\gamma$, and such that
$t_\beta=p_\beta$ and $t_\gamma=q_\gamma$. Assume moreover that
for $\alpha=\beta$ and for $\alpha=\gamma$,
$p_\alpha,\,q_\alpha,\,r_\alpha$ and $s_\alpha$ are all different
and that $p_\alpha\bigvee q_\alpha$ covers
$p_\alpha,\,q_\alpha,\,r_\alpha$ and $s_\alpha$. By Definition
\ref{DefinitionSeparatedProduct}, $a=\{p,q,r\}$ and
$b=\{p,q,r,s\}$ are in $\ovee_\alpha\l_\alpha$. Let $R^0=a\bigcup
t$. Then (see Lemma \ref{LemmaJoininTensorF}),
\[\begin{aligned}R^1&:=\bigvee_\gamma R^0=R^0\cup p[p_\gamma\vee q_\gamma]\\
R^2&:=\bigvee_\beta R^1=R^0\cup p[p_\gamma\vee
q_\gamma]\cup q[p_\beta\vee q_\beta]\cup r[p_\beta\vee
q_\beta] \\
R^3&:=\bigvee_\gamma R^2=p[p_\beta\vee q_\beta,p_\gamma\vee
q_\gamma]\, .\end{aligned}\]
Hence, by Lemma \ref{LemmaJoininTensorF}, $a\bigvee_\ovee t=R^3
\supsetneq b\supsetneq a$, therefore $\ovee_\alpha\l_\alpha$ does not have
the covering property.
\end{proof}

\begin{definition}
Let $\s_1,\,\s_2$ be sets, $\l_1=\mo_{\s_1}$, and
$\l_2=\mo_{\s_2}$ (see Definition \ref{DefinitionXi}). Then,
\[\l_1\circ\l_2:=\l_1\owedge\l_2\cup\{R\in\Xi(\s_1\t\s_2)
\ \lvert\ \abs R=3\}.\]
\end{definition}

\begin{theorem}
Let $\s_1$ and $\s_2$ be sets, $\l_1=\mo_{\s_1}$,
$\l_2=\mo_{\s_2}$, $T=\aut(\l_1)\t\aut(\l_2)$, and
$\l\in\Sep_T(\l_1,\l_2)$. Suppose that for $i=1$ and 2,
$\vert\s_i\vert=3$ or $\vert \s_i\vert\geq 5$. Then $\l$ has the covering property
if and only if $\l=\l_1\circ\l_2$.
\end{theorem}

\begin{proof} ($\Leftarrow$): Let $\mathbf{\s}=\s_1\t\s_2$,
$\Xi=\Xi(\s_1\t\s_2)$, and
$a\in\l_1\circ\l_2$. Then $a\in\Xi$ and $\abs{a}=3$, or
$a\in\l_1\owedge\l_2$. Hence one of the following cases holds.
\begin{enumerate}
\item $a\in\mathbf{\s}$.
\item $a\in\Xi$ and $\abs{a}=2$ or $3$.
\item $a=p_1\t\s_2$ or $a=\s_2\t p_2$ for some
$p=(p_1,p_2)\in\mathbf{\s}$.
\item $a=p_1\t\s_2\cup\s_1\t p_2$, for some $p\in\mathbf{\s}$
({\it i.e.}, $a$ is a coatom).
\end{enumerate}
Hence, obviously $\l_1\circ\l_2$ has the covering property.

($\Rightarrow$): Let $\l\in\sep(\l_1,\l_2)$ with the covering
property, and $R\subseteq\mathbf{\s}$.

(1): By Lemma \ref{LemmaJoininTensorF} and Axiom P2, if
$R\not\in\Xi$ ({\it i.e.}, there are $p^1,\,p^2\in R$ with
$p^1_i=p^2_i$ for $i=1$ or $2$), then $\bigvee
R\in\l_1\owedge\l_2$.

(2): Suppose now that $R\in\Xi$. By Lemma
\ref{LemmaJoinP1diffQ1andP2diffQ2} part 1, if $\abs R\leq 2$, then
$R\in\l_1\owedge\l_2$, hence $R\in\l$. Moreover, if $\abs R\geq 3$,
then for all $s\in\mathbf{\s}$ not in $R$ with $s_1\in \pi_1(R)$
or $s_2\in\pi_2(R)$, we have, by Lemma \ref{LemmaJoininTensorF}, $\bigvee (R\bigcup s)=1$.

(3) {\bf Claim}: Suppose that $R\in\Xi$ and that $\abs R=3$. Write $a:=\bigvee
R$. Then $a\ne 1$ 

{\it Proof}\nobreak : Write
$R=\{p,q,r\}$ and suppose that $a=1$. As we have seen
$R^0:=\{p,q\}\in\l$. Hence $1=a=r\vee R^0\supsetneq p_1\t\s_2\cup\s_1\t q_2\supsetneq R^0$,
a contradiction since $\l$ has the covering property. This proves the claim.

As a consequence, from (2), $1\gtrdot a$. Moreover, $a\in\Xi$.
We write $a$ as $a=\{p_1\t f(p_1)\ \lvert\ p_1\in\pi_1(a)\}$. Hence,
$f$ is injective.

(4) {\bf Claim}: $a=R$. 

{\it Proof}\nobreak : 
If $\vert\s_1\vert=3$ or $\vert\s_2\vert=3$, the proof is finished. 
So we can assume that $\vert\s_i\vert \geq 5$ for
$i=1,2$.
Note that
any bijection of $\s_2$ induces an automorphism of $\l_2$. Suppose
that $a\ne R$, hence $\abs{a}\geq 4$. Let $v_2\in\aut(\l_2)$ such
that its restriction to $\pi_2(a)$ is different from the identity,
and with at least three fixed points in $\pi_2(a)$.  By Axiom P4,
there is $u\in\aut(\l)$ such that on $\mathbf{\s}$, $u$ equals
$\mathsf{id}\t v_2$. Hence $c:=\{(p_1,v_2\circ f(p_1))\ \lvert\ 
p_1\in\pi_1(a)\}\in\l$ and $c\cap a\ne a$, therefore $1$ does
not cover $c\cap a$. Moreover $\abs{c\cap a}\geq 3$; whence
a contradiction by (3). This proves the claim.

(5) Finally, suppose that $R\in\Xi$ and $\abs{R}\geq 4$. Let
$b=\bigvee R$. To see that $b=1$, let
$R_0\in\Xi$ with $R_0\subseteq R$ and $\abs{R_0}=3$. By what
precedes, $R_0\in\l$ and $1\gtrdot R_0$. As a consequence, $b=1$.
\end{proof}
\section{An example with the covering
property}\label{SectionAnExample}

Let $\h_1$ and $\h_2$ be complex Hilbert spaces and
$T=\uni(\h_1)\t\uni(\h_2)$. In this section, we give an example in
$\Sep_T(\pro(\h_1),\pro(\h_2))$, denoted by
$\pro(\h_1)\ots\pro(\h_2)$, which has the covering property.
Moreover, $\pro(\h_1)\ots\pro(\h_2)$ is coatomistic, but, as
expected, the dual has not the covering property, {\it i.e.
$\pro(\h_1)\ots\pro(\h_2)$ is not a DAC-lattice.}

\begin{definition}\label{DefinitionOTS}
Let $\h_1$ and $\h_2$ be complex Hilbert spaces,
$\mathbf{\s}=\s_{\h_1}\t\s_{\h_2}$, $\l_1=\pro(\h_1)$,
$\l_2=\pro(\h_2)$, and $\s_\otimes$ the set of one-dimensional
subspaces of $\h_1\otimes\h_2$. For $V\in\pro(\h_1\otimes\h_2)$, define
$\prs[V]:=\{(p_1,p_2)\in\mathbf{\s}\ \lvert\ p_1\otimes p_2\in V\}$. Then, 
$\l_1\ots\l_2:=\{\prs[V]\ \lvert\  V\in \pro(\h_1\otimes \h_2)\}$,
ordered by set-inclusion. For $A\subseteq\mathbf{\s}$, we write
$A^\p:=\{q\in\s_\otimes\ \lvert\ \langle q,p_1\otimes p_2\rangle=0,\,\forall
(p_1,p_2)\in A\}$, where $\langle-\vert-\rangle$ denotes the
scalar product in $\h_1\otimes\h_2$. Moreover, we denote the set of
antilinear maps from $\h_1$ to $\h_2$ by $\mathcal{A}(\h_1,\h_2)$.
\end{definition}

\begin{proposition} Let $m$ and $n$ be
integers, $\h_1=\mathbb{C}^m$, $\h_2=\mathbb{C}^n$,
$\l_1=\pro(\h_1)$, $\l_2=\pro(\h_2)$, and
$\mathbf{\s}=\s_{\h_1}\t\s_{\h_2}$. For $A\in
\mathcal{A}(\h_1,\h_2)$, define $X_A\subseteq\mathbf{\s}$ as
$X_A:=\bigcup\{p_1\t(A(p_1)^{\p_2})\ \lvert\ p_1\in \s_1\}$. Then,
\[\l_1\ots\l_2=\{\bigcap\omega\ \lvert\ \omega\subseteq\{X_A\ \lvert\ A\in
\mathcal{A}(\h_1,\h_2)\}\}.\]
\end{proposition}

\begin{proof}
Let $\{e_i^1\}_{1\leq i\leq m}$ and $\{e_j^2\}_{1\leq j\leq n}$
denote the canonical basis of $\mathbb{C}^m$ and $\mathbb{C}^n$
respectively.

Let $v\in\mathbb{C}^m\otimes\mathbb{C}^n$ and
$p=(p_1,p_2)\in\mathbf{\s}$ with $p_1=\mathbb{C}w_1$ and
$p_2=\mathbb{C}w_2$. Write $v$, $w_1$ and $w_2$ as
\[v=\sum_{i=1}^m\sum_{j=1}^n s_{ij} e_i^1\otimes e_j^2,\
w_1=\sum_{i=1}^m\lambda_i e_i^1,\ \mbox{and}\ w_2=\sum_{j=1}^n
\mu_j e_j^2.\]
Let $\lambda=(\lambda_1,\cdots,\lambda_m)^T\ \mbox{and}\
\mu=(\mu_1,\cdots,\mu_n)^T$.
Let $S$ be the $m\t n$ matrix defined as $S_{ij}=s_{ij}$. Then
$p\in \prs[v^\p]$ if and only if $\langle w_1\otimes w_2,v\rangle=0$,
hence if and only if $\overline{\mu}^T S^T
\overline{\lambda}=0$. Let $A$ be the antilinear map defined by
the matrix $S^\p$. Then, $p\in \prs[v^\p]$ if and only if
$p_2\in A(p_1)^\p$, that is if and only if $p\in X_A$. As a
consequence, $\prs[v^\p]=X_A$.

On the other hand, if $A\in \mathcal{A}(\h_1,\h_2)$, then
$X_A=\prs[v^\p]$, where $v$ is given by the formula above with
$s_{ij}=(A^T)_{ij}$.
\end{proof}

\begin{remark}
Let $\h$ be a complex Hilbert space of dimension $\geq 3$. Then,
by Wigner's theorem (see \cite{Faure/Froelicher:handbook}, Theorem
14.3.6), any ortho-automorphism of $\pro(\h)$ is induced by a
unitary or antiunitary map on $\h$. Note that if $v_1$ is a
unitary map on $\h_1$ and $v_2$ is an antiunitary map on $\h_2$,
then $v=v_1\t v_2$ does not induce an automorphism of
$\l_1\ots\l_2$. Indeed, let $X_A$ be a coatom. Then
\[v(X_A)=\bigcup_{p_1\in\s_1}v_1(p_1)\t
v_2(A(p_1)^{\p_2})=\bigcup_{p_1\in\s_1}p_1\t((v_2\circ A\circ
v_1^{-1}(p_1))^{\p_2}).\]
Now, since $A$ and $v_2$ are antilinear and $v_1$ is linear, it
follows that $v_2\circ A\circ v_1^{-1}$ is linear, hence $v(X_A)$
is not a coatom of $\l_1\ots\l_2$.
\end{remark}

\begin{theorem}\label{TheoremP(H1)AertsP(H2)}
Let $\h_1$ and $\h_2$ be complex Hilbert spaces,
$\mathbf{\s}=\s_{\h_1}\t\s_{\h_2}$, $\l_1=\pro(\h_1)$,
$\l_2=\pro(\h_2)$, and $T=\uni(\h_1)\t\uni(\h_2)$. Then
\begin{enumerate}
\item for all $A\subseteq\mathbf{\s}$, we have $\bigvee_\otssub
A=\prs[A^{\p\p}]$,
\item $\l_1\ots\l_2\in\Sep_T(\l_1,\l_2)$,
\item $\l_1\ots\l_2$ has the covering property and is
coato\-mistic, but, if $\l_1\ne 2\ne\l_2$ ({\it i.e.}, the
dimension of $\h_1$ and $\h_2$ is $\geq 2$), the dual has not the
covering property,
\item $\l_1\owedge\l_2=\{\prs[V]\ \lvert\  V\in\pro(\h_1\otimes \h_2),\,
V=\prs[V]^{\p\p}, V^\p=\prs[V^\p]^{\p\p}\}$ ({\it i.e.}, both $V$
and $V^\p$ are spanned by product vectors),
\item if $\l_1\ne 2\ne\l_2$, then $\l_1\owedge\l_2\varsubsetneqq
\l_1\ots\l_2\varsubsetneqq\l_1\ovee\l_2$.
\end{enumerate}
\end{theorem}

\begin{proof} (1): This follows directly from Definition \ref{DefinitionOTS}.

(2): Obviously, $\l_1\ots\l_2$ is a simple closure space on $\mathbf{\s}$. Let
$a\in\l_1\owedge\l_2$. By Definition, $a^\#=\prs[a^\p]$.
Hence, $a^\#\subseteq a^\p$, $a^{\p\p}\subseteq a^{\#\p}$, thus
$\prs[a^{\p\p}]\subseteq a^{\#\#}=a$, and therefore
$a=\prs[a^{\p\p}]$. As a consequence,
$\l_1\owedge\l_2\subseteq\l_1\ots\l_2$.

Let $V\in\pro(\h_1\otimes\h_2)$ such that $\prs[V]=p_1\t A_2$. Then
\[p_1\otimes(A_2^{\p_2\p_2})=(p_1\otimes A_2)^{\p\p}\subseteq V,\]
therefore $A_2^{\p_2\p_2}\subseteq A_2$, hence $A_2\in\l_2$. As a
consequence, Axiom P3 holds.

Axiom P4 with $T=\uni(\h_1)\t\uni(\h_2)$ holds in
$\pro(\h_1\otimes\h_2)$, therefore obviously also in $\l_1\ots\l_2$.

(3): The covering property holds in $\pro(\h_1\otimes \h_2)$ (see
\cite{Maeda/Maeda:handbook}, Theorem 34.2), hence, by (1), also
in $\l_1\ots\l_2$. Moreover, since $\pro(\h_1\otimes\h_2)$ is
coatomistic, so is $\l_1\ots\l_2$.

Next, let $p\in\mathbf{\s}$. Then $x=\prs[p^\p]=p^{\#}$ is a
coatom of $\l_1\ots \l_2$. Now, there is $R\in\Xi(\mathbf{\s})$
(see Definition \ref{DefinitionXi}) with $\abs R=2$, such that
$x\cap R=\emptyset$. By Lemma
\ref{LemmaJoinP1diffQ1andP2diffQ2}, $R\in\l_1\ots\l_2$. Moreover
$x\bigvee_\otssub R=1$ since $x$ is a coatom. Hence, writing
$R=\{p,q\}$, and the order relation, meet, join, bottom and top
elements in $(\l_1\ots\l_2)^*$ by $\leq_*$, $\bigwedge_*$,
$\bigvee_*$, $0_*$, and $1_*$ respectively, we find that
$x\bigwedge_* R=0_*$ and $x\bigvee_* R=1_*\gneqq_* p\gneqq_* R$.
Therefore, $(\l_1\ots\l_2)^*$ does not have the covering property.

(4) Let $p\in\mathbf{\s}$ and $q\in p^{\#\p}$. Write
$q=\mathbb{C}v$ with $v\in\h_1\otimes\h_2$. For $i=1$ and $i=2$, let
$\{w_i^k\}$ be an ortho-basis of $p_i^{\p_i}$, and let $x_i\in
p_i$ ({\it i.e.}, $p_i=\mathbb{C}x_i$). Then $v$ can be decomposed
as
\[v=\alpha x_1\otimes x_2+\sum_{k_2}\beta_{k_2}x_1\otimes
w_2^{k_2}+\sum_{k_1}\beta_{k_1}w_1^{k_1}\otimes
x_2+\sum_{l_1l_2}\gamma_{l_1l_2}w_1^{l_1}\otimes w_2^{l_2}.\]
Now, $p^\#=p_1^{\p_1}\t\s_{\h_2}\cup\s_{\h_1}\t p_2^{\p_2}$.
Hence, since $q\in p^{\#\p}$, we find that
$\gamma_{l_1l_2}=\beta_{k_1}=\beta_{k_2}=0$, for all
$k_1,\,k_2,\,l_1$ and $l_2$. Therefore, $v\in p$, hence
$p^{\#\p}=p$.

Let $a\in\l_1\owedge\l_2$. From (1), $a=\prs[a^{\p\p}]$. On the
other hand, $a^\#=\prs[a^\p]$. Now, $a=a^{\#\#}=\bigcap\{p^\#\ \lvert\ 
p\in a^\#\}$. Hence, by the preceeding,
\[a^\p=\lp \bigcup\{ p^{\#\p}\ \lvert\  p\in a^\#\}\rp^{\p\p}=
a^{\#\p\p},\]
therefore $a^\p$ is also spanned by product vectors. Hence, writing $V=a^{\p\p}$,
we find that $a=\prs[V]$, $V=\prs[V]^{\p\p}$ and $V^\p=\prs[V^\p]^{\p\p}$.

Let $V\in\pro(\h_1\otimes \h_2)$ such that both $V$ and $V^\p$ are
spanned by product vectors. Let $a:=\prs[V]$ and $b:=\prs[V^\p]$.
Since $\prs[V]^{\p\p}=V$ and also $\prs[V^\p]^{\p\p}=V^\p$, then
$\prs[\prs[V]^\p]= \prs[V^\p]\  \mbox{and}\ \prs[\prs[V^\p]
^\p]=\prs[V]$.
Therefore $a^\#=b$ and $b^\#=a$. As a consequence,
$a\in\l_1\owedge\l_2$.

(5): By (4), $\l_1\owedge\l_2\ne\l_1\ots\l_2$. On the other hand,
by Theorem \ref{TheoremCoveringinTensorF}, $\l_1\ovee\l_2$ does not have
the covering property, whereas by (3), $\l_1\ots\l_2$ has the
covering property. As a consequence,
$\l_1\ots\l_2\ne\l_1\ovee\l_2$.\end{proof}

\end{document}